\documentclass[11pt,a4paper,equation]{article}
\usepackage[greek,american]{babel}
\usepackage[iso-8859-7]{inputenc}
\usepackage{amssymb,amsfonts}
\usepackage[dvips]{graphicx}
\usepackage{amsmath}
\usepackage{epsfig}
\usepackage{latexsym}
\usepackage{makeidx}
\usepackage{pst-math,pst-xkey}

\numberwithin{equation}{section}
\makeindex
\begin{document}
\date{}
\author{\textbf{Vassilis G. Papanicolaou$^1$ and Aristides V. Doumas$^2$}
\\\\
Department of Mathematics
\\
National Technical University of Athens,
\\
Zografou Campus, 157 80, Athens, GREECE
\\\\
$^1${\tt papanico@math.ntua.gr} \qquad  $^2${\tt adou@math.ntua.gr}}
%^{2}$aris.doumas@hotmail.com \quad
\title{The Coupon Collector's Brother}
\maketitle
\begin{abstract}
A popular variant of the collector's problem is the following:
Assume there are $N$ different types of
coupons with equal occurring probabilities. 
There is one main collector who collects coupons.
When she gets a ``double," she gives it to her younger brother. Hence,
when the main collector completes her collection, her brother's album
will still have, say, $U_N$ empty spaces. We show that, as $N \to \infty$, the limiting distribution of $U_N/\ln N$ is exponential
with parameter $1$.
\end{abstract}
\textbf{2010 AMS Mathematics Classification.} 60F05; 60F99.
\section{Introduction}

The classical \emph{coupon collector's problem (CCP)} concerns a population
(e.g. fishes, viruses, genes, words, baseball cards, etc.)
whose members are of $N$ different types. The members of the population are
sampled independently with replacement and their types are recorded. CCP pertains to the family of urn problems along with other famous problems, such as the birthday problem and the occupancy problem. Its origin can be traced back to De Moivre's treatise \textit{De Mensura Sortis} of 1712 (see,
e.g., \cite{Ho}) and Laplace's pioneering work \textit{Theorie Analytique de Probabilites} of 1812 (see \cite{D-H}).
CCP has attracted the attention of various researchers due to the fact that it has applications
%(and implications)
in many areas of science (computer science/search algorithms, mathematical programming, optimization, learning processes, engineering, ecology, as well as linguistics---see, e.g., \cite{B-H}).\\

Some preliminary results regarding the CCP were presented in the two relatively unknown works \cite{M} and \cite{G} (where the entertaining term \emph{cartophily} appeared in the titles). However, the  classical references of the CCP are considered to be
the well-known book \cite{F} of Feller (Chapter IX, pp. 225, 239), the paper \cite{N-S} of Newman and Shepp on the \emph{Double Dixie Cup
problem}, where the authors compute the asymptotics, as $N \to \infty$, of the expected value $\mathbb{E}[T_N]$,
$T_N$ being is the number of trials it takes in order to obtain $m$ complete sets of $N$ coupons, and the paper
\cite{E} of Erd\H{o}s and R$\acute{e}$nyi, where the limit distribution of the random variable $T_N$ (appropriately normalized) has been derived.

An interesting version of the classical CCP can be described as follows:\footnote{For an entertaining computer simulation see \cite{PS}.}
Assume that a
company sells ice cream with a cardboard cover that has hidden on the underside a picture (``coupon'') of a sixties music band. In total there are $N$ different pictures and each one appears with probability $1/N$. Mr. and Mrs. Smith have one daughter and one son, both ice cream and sixties music addicts. The girl (she is the oldest) is the only one to buy ice cream and tries to complete her collection. When she gets a new picture she puts it
in her album, and when she gets a double she gives it to her brother. After having bought $T_N$ ice creams, the girl has completed her album, while they remain $U_N$ unfilled places in her brother's album, where, obviously, $1 \leq U_N \leq N$. The problem here is to study the random variable
$U_N$. And one major aspect of this problem is the behavior of $U_N$ as $N$ becomes large.

This CCP variant has been studied for quite a while. Pintacuda \cite{P} (see also \cite{Ro}), by using the martingale stopping theorem, proved that
\begin{equation}
\mathbb{E}\left[U_N\right] = \sum_{m=1}^N \frac{1}{m} =: H_N
\label{00a}
\end{equation}
($H_N$ is sometimes called the $N$-th harmonic number). Thus,
\begin{equation}
\mathbb{E}[U_N] \sim \ln N
\qquad \text{as }\; N \to \infty.
\label{00b}
\end{equation}
Foata, Guo-Niu, and Lass  \cite{FO} and Foata and Zeilberger \cite{FD}, using nonelementary mathematics, extended \eqref{00a}
to the case of many brothers. Soonafter, Adler, Oren, and Ross \cite{Ad} derived the same results by using basic probability arguments.

By exploiting the techniques of \cite{Ad} one can compute explicitly the variance of $U_N$ as
\begin{equation}
\mathbb{V}[U_N] = 4 \left(H_1 + \frac{H_2}{2} + \cdots + \frac{H_N}{N}\right) - 3 H_N - H_N^2
\label{A0a}
\end{equation}
and from \eqref{A0a} it follows easily that
\begin{equation}
\mathbb{V}[U_N] \sim \ln^2 N
\qquad \text{as }\; N \to \infty.
\label{A0b}
\end{equation}

A question that comes naturally from the above formulas is \cite{D-P} what is the limiting distribution of $U_N$ as $N \to \infty$. The nature of the random variable $U_N$ together with the asymptotic formulas \eqref{00b} and \eqref{A0b} hint that 
\begin{equation}
\frac{U_N}{\ln N}   \, \overset{d}{\longrightarrow} \, \mathcal{E}(1)
\qquad \text{as }\; N \to \infty,
\label{A0}
\end{equation}
where $\mathcal{E}(1)$ is an exponential random variable with parameter $1$, while, as usual, $\overset{d}{\longrightarrow}$ denotes convergence in
distribution.

The purpose of this note is to show that \eqref{A0} is, actually, true.

\section{Proof of \eqref{A0}}
Before proceeding with the proof, let us recall the well-known formula
\begin{equation}
e^{-x} \sum_{\ell = 0}^{L-1} \frac{x^{\ell}}{\ell !} = \frac{1}{(L-1)!}\int_x^{\infty} y^{L-1} e^{-y} dy.
\label{A2a}
\end{equation}
For $x > 0$ formula \eqref{A2a} reflects the relation between the Poisson and the Erlang distributions (while another immediate way to check its validity
is by differentiating both sides with respect to $x$). Now, from \eqref{A2a} we get
\begin{equation}
\sum_{\ell = 0}^{L-1} \frac{x^{\ell}}{\ell !} = \frac{1}{(L-1)!}\int_x^{\infty} y^{L-1} e^{-(y-x)} dy
\label{A2}
\end{equation}
and, finally, by substituting $\xi = y-x$ in the integral of \eqref{A2} we obtain
\begin{equation}
\sum_{\ell = 0}^{L-1} \frac{x^{\ell}}{\ell !} = \frac{1}{(L-1)!}\int_0^{\infty} (x + \xi)^{L-1} e^{-\xi} d\xi.
\label{A1}
\end{equation}
Actually, by analytic continuation it follows automatically that \eqref{A1} is valid for any complex number $x$.

\medskip

\textit{Proof of formula \eqref{A0}}. Our starting point is the formula
\begin{equation}
G_N(s) := \mathbb{E}\left[s^{U_N}\right]
= s \sum_{k+\ell+m = N-1} \binom{N-1}{k,\ell,m} (-1)^{\ell} \frac{N(s-1)^m}{(N-k)^{m+1}} \, m!,
\label{A3}
\end{equation}
which has been derived by Foata, Guo-Niu, and Lass \cite{FO}. Here the sum is taken over all nonnegative integers $k$, $\ell$, and $m$ satisfying
$k+\ell+m = N-1$ and $\binom{N-1}{k,\ell,m}$ stands, as usual, for the trinomial coefficient, thus
\begin{equation*}
\binom{N-1}{k,\ell,m} = \frac{(N-1)!}{k! \, \ell ! \, m!}.
%\label{A3a}
\end{equation*}
Clearly, $G_N(s)$ is a polynomial of degree $N$.

Formula \eqref{A3} implies (since $m = N - 1 - k - \ell$)
\begin{align}
\frac{G_N(s)}{s}
&= N! \sum_{k=0}^{N-1} \sum_{\ell=0}^{N-1-k} \frac{(-1)^{\ell}}{k! \, \ell !} \frac{(s-1)^{N-1-k-\ell}}{(N-k)^{N-k-\ell}}
\nonumber
\\
&= N! \sum_{k=0}^{N-1} \frac{1}{k!} \frac{(s-1)^{N-1-k}}{(N-k)^{N-k}}
\left[\sum_{\ell=0}^{N-1-k} \frac{1}{\ell !} \frac{(N-k)^{\ell} }{(1-s)^{\ell} }\right]
\label{A4}
\end{align}
Next, we apply formula \eqref{A1} to \eqref{A4} (for $L = N-k$ and $x = \frac{N-k}{1-s}$) and obtain
\begin{align}
\frac{G_N(s)}{s}
&= N! \sum_{k=0}^{N-1} \frac{1}{k!} \frac{(s-1)^{N-1-k}}{(N-k)^{N-k}}
\frac{1}{(N-k-1)!} \int_0^{\infty} \left(\frac{N-k}{1-s} + \xi\right)^{N-k-1} e^{-\xi} d\xi
\nonumber
\\
&= \sum_{k=0}^{N-1} \binom{N}{k} \frac{(s-1)^{N-1-k}}{(N-k)^{N-k-1}}
\int_0^{\infty} \left(\frac{N-k}{1-s} + \xi\right)^{N-k-1} e^{-\xi} d\xi
\nonumber
\\
&= \int_0^{\infty} \left\{\sum_{k=0}^{N-1} \binom{N}{k} \frac{(-1)^{N-1-k}}{(N-k)^{N-k-1}}
\left[N-k + (1-s)\xi\right]^{N-k-1}\right\} e^{-\xi} d\xi.
\label{A5}
\end{align}
By replacing the summation (dummy) index $k$ by $N-k$ in the above sum we obtain
\begin{align}
\frac{G_N(s)}{s}
&= \int_0^{\infty} \left\{\sum_{k=1}^N \binom{N}{k} \frac{(-1)^{k-1}}{k^{k-1}}
\left[k + (1-s)\xi\right]^{k-1}\right\} e^{-\xi} d\xi
\nonumber
\\
&= \int_0^{\infty} \left\{\sum_{k=1}^N (-1)^{k-1} \binom{N}{k}
\left[1 + \frac{(1-s)\xi}{k}\right]^{k-1}\right\} e^{-\xi} d\xi.
\label{A6}
\end{align}
If there were not a $k$ appearing inside the square brackets in the last expression of \eqref{A6}, then we could easily compute the sum. Fortunately,
we can get rid of the $k$ by substituting $\xi = kx$. Then, \eqref{A6} yields
\begin{equation}
\frac{G_N(s)}{s}
= \int_0^{\infty} \left\{\sum_{k=1}^N (-1)^{k-1} k \binom{N}{k}
\left[1 + (1-s) x\right]^{k-1} e^{-(k-1) x}\right\} e^{-x} dx
\label{A7}
\end{equation}
and it is easy to see that
\begin{equation*}
\sum_{k=1}^N (-1)^{k-1} k \binom{N}{k} \left[1 + (1-s) x\right]^{k-1} e^{-(k-1) x}
= N \left\{1 - \left[1 + (1-s) x\right] e^{-x}\right\}^{N-1}.
%\label{A8}
\end{equation*}
Therefore, formula \eqref{A7} becomes
\begin{equation}
\frac{G_N(s)}{s}
= N \int_0^{\infty} \left\{1 - \left[1 + (1-s) x\right] e^{-x}\right\}^{N-1} e^{-x} dx
\label{A9}
\end{equation}
or, by integrating by parts,
\begin{equation}
\frac{G_N(s)}{s}
= (1-s) \int_0^{\infty} \left\{1 - \left[1 + (1-s) x\right] e^{-x}\right\}^N \frac{dx}{\left[s + (1-s) x\right]^2}.
\label{A10}
\end{equation}
Notice that formula \eqref{A10} is valid for any $s \in \mathbb{C}$ ($s=1$ is a removable singularity of the right-hand side, while, obviously,
$G_N(1) = 1$).

Since we are interested in the limiting distribution of the random variable $U_N / \ln N$, we can consider its characteristic function $\phi_N(t)$
and try to find its limit (for any fixed $t \in \mathbb{R}$) as $N \to \infty$. In view of \eqref{A3} we have
\begin{equation}
\phi_N(t) = \mathbb{E}\left[e^{i t \, U_N / \ln N}\right] = G_N\left(e^{i t / \ln N}\right)
\label{A11}
\end{equation}
and, then, by \eqref{A10} and \eqref{A11} we get
\begin{equation}
\frac{\phi_N(t)}{e^{i t / \ln N}}
= \left(1 - e^{i t / \ln N}\right) \int_0^{\infty}
\frac{\left\{1 - \left[1 + \left(1 - e^{i t / \ln N}\right) x\right] e^{-x}\right\}^N}
{\left[e^{i t / \ln N} + \left(1 - e^{i t / \ln N}\right) x\right]^2} \, dx.
\label{A12}
\end{equation}
Now, the substitution $x = \xi \ln N$ in the integral of \eqref{A12} yields
\begin{equation}
\frac{\phi_N(t)}{e^{i t / \ln N}}
= \left(1 - e^{i t / \ln N}\right) \big[I_1(N) + I_2(N) + I_3(N)\big] \ln N,
\label{A13}
\end{equation}
where we have set
\begin{equation}
I_1(N) := \int_0^{1-\varepsilon}
\frac{\left\{1 - \left[1 + \xi \left(1 - e^{i t / \ln N}\right) \ln N \right] N^{-\xi}\right\}^N}
{\left[e^{i t / \ln N} + \xi \left(1 - e^{i t / \ln N}\right) \ln N \right]^2} \, d\xi,
\label{A13a}
\end{equation}
\begin{equation}
I_2(N) := \int_{1-\varepsilon}^{1+\varepsilon}
\frac{\left\{1 - \left[1 + \xi \left(1 - e^{i t / \ln N}\right) \ln N \right] N^{-\xi}\right\}^N}
{\left[e^{i t / \ln N} + \xi \left(1 - e^{i t / \ln N}\right) \ln N\right]^2} \, d\xi,
\label{A13b}
\end{equation}
and
\begin{equation}
I_3(N) := \int_{1+\varepsilon}^{\infty}
\frac{\left\{1 - \left[1 + \xi \left(1 - e^{i t / \ln N}\right) \ln N \right] N^{-\xi}\right\}^N}
{\left[e^{i t / \ln N} + \xi \left(1 - e^{i t / \ln N}\right) \ln N\right]^2} \, d\xi.
\label{A13c}
\end{equation}
Here $\varepsilon$ is any number in $(0, 1)$.

Let us, first, notice that
\begin{equation}
\left(1 - e^{i t / \ln N}\right) \ln N = -it + O\left(\frac{1}{\ln N}\right),
\qquad
N \to \infty.
\label{A14}
\end{equation}
Thus, by using \eqref{A14} in \eqref{A13} we get immediately that
\begin{equation}
\lim_N \phi_N(t) = -it \left[\lim_N I_1(N) + \lim_N I_2(N) + \lim_N I_3(N)\right],
\label{A15}
\end{equation}
provided all limits exist.

Next, let us consider the numerator of the integrands in \eqref{A13a}, \eqref{A13b}, and \eqref{A13c}, namely the expression
\begin{equation}
\left[1 - \frac{1 + \xi \left(1 - e^{i t / \ln N}\right) \ln N }{N^{\xi}} \right]^N
=e^{A_N(\xi)},
\label{A16a}
\end{equation}
where for typographical convenience we have set
\begin{equation}
A_N(\xi) := N \ln \left(1 - \frac{1 + \xi \left(1 - e^{i t / \ln N}\right) \ln N }{N^{\xi}}\right).
\label{A16b}
\end{equation}
In view of \eqref{A14}, formula \eqref{A16b} yields
\begin{equation}
A_N(\xi) = N \ln \left(1 - \frac{1 - it\xi + O\left(\ln^{-1} N\right) \xi}{N^{\xi}}\right),
\qquad
N \to \infty.
\label{A17}
\end{equation}
(i) Suppose $0 \leq \xi \leq 1-\varepsilon$. Then, \eqref{A17} implies
\begin{equation}
A_N(\xi) = -N^{1-\xi} + it\xi N^{1-\xi} + O\left(\frac{N^{1-\xi}}{\ln N}\right),
\qquad
N \to \infty,
\label{A18}
\end{equation}
uniformly in $\xi$. Consequently,
\begin{equation}
\lim_N e^{A_N(\xi)} = 0,
\qquad \text{uniformly in }\; [0, 1-\varepsilon].
\label{A19}
\end{equation}
Since the denominator of the integrand in \eqref{A13a} is bounded away from $0$ for all sufficiently large $N$, we can conclude from \eqref{A16a}
and \eqref{A19} that
\begin{equation}
\lim_N I_1(N) = 0.
\label{A20}
\end{equation}

(ii) Next, suppose $\xi \geq 1+\varepsilon$. Then, \eqref{A17} implies
\begin{equation}
A_N(\xi) = -\frac{1}{N^{\xi-1}} + \frac{it\xi}{N^{\xi-1}} + O\left(\frac{\xi}{N^{\xi-1}\ln N}\right),
\qquad
N \to \infty.
\label{A21}
\end{equation}
Therefore,
\begin{equation}
\lim_N A_N(\xi) = 0,
\qquad \text{uniformly in }\; [1+\varepsilon, \infty)
\label{A22}
\end{equation}
and, consequently,
\begin{equation}
\lim_N e^{A_N(\xi)} = 1,
\qquad \text{uniformly in }\; [1+\varepsilon, \infty).
\label{A23}
\end{equation}
Then, by dominated convergence, and in view of \eqref{A16a} and \eqref{A23}, formula \eqref{A13c} gives
\begin{equation}
\lim_N I_3(N) = \int_{1+\varepsilon}^{\infty} \frac{1}{\left(1 -it \xi\right)^2} \, d\xi
= -\frac{1}{it}\cdot \frac{1}{1 - it(1+\varepsilon)}.
\label{A24}
\end{equation}

(iii) Finally, regarding the integral $I_2(N)$, we can observe that the integrand is bounded as $N \to \infty$. Therefore,
\begin{equation}
\left|I_2(N)\right| \leq M \varepsilon
\label{A25}
\end{equation}
where $M$ does not depend on $N$.

By using \eqref{A20}, \eqref{A24}, and \eqref{A25} in \eqref{A15}, and recalling that $\varepsilon \in (0, 1)$ is arbitrary, we can conclude that
\begin{equation}
\lim_N \phi_N(t) = \frac{1}{1 - it}
\qquad \text{for any }\; t \in \mathbb{R}
\label{A26}
\end{equation}
and the proof of \eqref{A0} is completed, since the right-hand side of \eqref{A26} is the characteristic function of the exponential distribution
with parameter $1$.
\hfill $\blacksquare$

%\medskip

% \textbf{Acknowledgments.}

\end{document}